\documentclass[12pt]{amsart}

\usepackage{latexsym}
\usepackage{amsfonts}
\usepackage{amssymb}
\usepackage{fullpage}
\usepackage[OT4]{fontenc}
\usepackage{amscd}
\usepackage{textcomp}
\usepackage{color}

\newtheoremstyle{s2}{9pt}{9pt}{\rm}{\parindent}{\bf}{.}{0.5em}{}
\theoremstyle{s2} 

\newtheorem{definition}{Definition}

\newtheoremstyle{s1}{9pt}{9pt}{\it}{\parindent}{\bf}{.}{0.5em}{}
\theoremstyle{s1}
\newtheorem{lemma}[definition]{Lemma}
\newtheorem{theorem}[definition]{Theorem}

\newtheorem{proposition}[definition]{Proposition}
\newtheorem{remark}[definition]{Remark}

\newtheorem*{conjecture*}{Conjecture}

\def\P2{{\mathbb{P}^2}}
\def\P1{{\mathbb{P}^1}}
\def\Ce{\mathbb{C}}

\def\Zet{\mathbb{Z}}
\def\En{\mathbb{N}}

\def\Qu{\mathbb{Q}}

\def\Oo{\mathcal{O}}

\DeclareMathOperator{\mult}{mult}
\DeclareMathOperator{\lcm}{lcm}
\DeclareMathOperator{\Num}{Num}

\numberwithin{definition}{section}

\title[A note on $k$-very ampleness]
{A note on $k$-very ampleness of line bundles on general blow-ups of hyperelliptic surfaces} \makeatletter

\@addtoreset{equation}{section}
\makeatother
\author{{\L}ucja Farnik}
\address{{\L}ucja Farnik, Jagiellonian University, Faculty of Mathematics and Computer Science, {\L}ojasiewicza~6, 30-348 Krak\'{o}w, Poland} 
\email{lucja.farnik@uj.edu.pl}

\date{\today}

\begin{document}
\bibliographystyle{alpha}

\begin{abstract}
We study $k$-very ampleness of line bundles on blow-ups of hyperelliptic surfaces at $r$ very general points. We obtain a numerical condition on the number of points for which a line bundle on the blow-up of a hyperelliptic surface at these $r$ points gives an embedding of order $k$.
\end{abstract}

\maketitle

\thispagestyle{empty}
\section{Introduction}

M.C. Beltrametti, P. Francia and A.J. Sommese introduced and studied
the concepts of higher order embeddings: $k$-spandness, $k$-very ampleness and $k$-jet ampleness of polarised varieties in a series of papers, see \cite{BeFS1989}, \cite{BeS1988}, \cite{BeS1993}. 
The problem of $k$-very ampleness on certain surfaces was studied by many authors. M. Mella and M. Palleschi in \cite{MP1993} proved the necessary and sufficient condition for a line bundle on any hyperelliptic surface to be $k$-very ample. Such a condition for any Del Pezzo surface was given by S.~Di Rocco in \cite{DR1996}.
Th. Bauer and T. Szemberg in \cite{BaSz1997} provided a criterion for $k$-very ampleness of a line bundle on an abelian surface.

In \cite{SzT-G2002} T. Szemberg and H. Tutaj-Gasi{\'n}ska  established a condition on the number of points for which a line bundle is $k$-very ample on a general blow-up of the projective plane. H. Tutaj-Gasi{\'n}ska  in \cite{T-G2002} gave a condition for $k$-very ampleness of a line bundle on a general blow-up of an abelian surface, and in \cite{T-G2005} --- on general blow-ups of elliptic quasi-bundles.

Recently, W. Alagal and A. Maciocia in \cite{AMa2014} study critical $k$-very ampleness on abelian surfaces, i.e. consider the critical value of
$k$ for which a line bundle is $k$-very ample but not $(k+1)$-very ample.

We come back to the classical question on the number of points for which a line bundle on a general blow-up of a surface is $k$-very ample. We consider blow-ups of hyperelliptic surfaces as such case has not been an object of study before.

\section{Notation and auxiliary results}
Let us set up the notation and basic definitions. We work over the field of complex numbers~$\Ce$. We consider only smooth reduced and irreducible projective varieties. By $D_1\equiv D_2$ we denote the numerical equivalence of divisors $D_1$ and $D_2$. By a curve we understand an irreducible subvariety of dimension 1. 
In the notation we follow \cite{Laz2004}.

\vskip 5pt
 We recall the definition of $k$-very ampleness.

Let  $X$ be a smooth projective variety of dimension $n$. Let  $L$ be a line bundle on  $X$, and let $x\in X$. 

\begin{definition}
We say that a line bundle $L$ is $k$-very ample if for every $0$-dimensional subscheme $Z\subset X$ of length $k+1$ the restriction map
 $$H^0(X,L)\longrightarrow H^0(X,L\otimes\Oo_Z)$$ is surjective.
\end{definition}

In the other words $k$-very ampleness means that the subschemes of length at most $k+1$ impose independent conditions on global sections of $L$.

We also recall the definition of the multi-point Seshadri constant.

Let  $x_1$, $\ldots$, $x_r\in X$ be pairwise distinct points.
\begin{definition}
The multi-point Seshadri constant of $L$ at $x_1$, $\ldots$, $x_r$
is the real number
$$\varepsilon(L,x_1, \ldots, x_r)=\inf \left\{\frac{LC}{\sum_{i=1}^r\mult_{x_i}C}: \ \{x_1,\ldots,x_r\}\cap C\neq \emptyset\right\},$$
where the infimum is taken over all irreducible curves $C\subset X$ passing through at least one of the points 
 $x_1$, $\ldots$, $x_r$.
\end{definition}

If $\pi\colon\widetilde{X}\longrightarrow X$ is the blow-up of  $X$ at  $x_1$, $\ldots$, $x_r$, and  $E_1$, $\ldots$, $E_r$ are exceptional divisors of the blow-up, then equivalently the Seshadri constant may be defined as (see e.g. \cite{Laz2004} vol. I, Proposition 5.1.5):
$$\varepsilon(L,x_1, \ldots, x_r)=\sup \left\{\varepsilon: \ \pi^*L-\varepsilon\sum_{i=1}^r E_i \text{ is nef} \right\}.$$

\vskip 5pt

Now let us recall the definition of a hyperelliptic surface.
\begin{definition}
A hyperelliptic surface $S$ (sometimes called bielliptic) 
is a surface with Kodaira  dimension equal to $0$ and  irregularity  $q(S)=1$.
\end{definition}
 
Alternatively (\cite{Bea1996}, Definition VI.19), a surface $S$ is hyperelliptic if  $S\cong (A\times B)/G$, where $A$ and $B$ are elliptic curves, and $G$ is an abelian group acting on A by translation and acting on B, such that  $A/G$ is an elliptic curve and $B/G\cong \mathbb{P}^1$; $G$ acts on $A\times B$ coordinatewise.
Hence we have the following situation:
$$
\begin{CD}
S\cong (A\times B)/G @>\Phi>> A/G @.\\
@V\Psi VV @.\\
B/G\cong \P1
\end{CD}
$$
where $\Phi$ and $\Psi$ are natural projections.

Hyperelliptic surfaces were classified at the beginning of 20th century by G. Bagnera and M.~de Franchis in \cite{BF1907}, and independently by F. Enriques i F. Severi in \cite{ES1909-10}. They showed that there are seven non-isomorphic types 
of hyperelliptic surfaces. Those types are characterised by the action of $G$ on $B\cong\Ce/(\Zet\omega\oplus\Zet)$ (for details see e.g. \cite{Bea1996}, VI.20). The canonical divisor $K_S$ of any hyperelliptic surface is numerically trivial.

In 1990 F.~Serrano in \cite{Se1990}, Theorem 1.4, characterised the group of classes of numerically equivalent divisors $\Num(S)$
for each of the  surface's type:
\begin{theorem}[Serrano] A basis of the group  $\Num(S)$
for each of the hyperelliptic surface's type and the multiplicities of the singular fibres in each case are the following:
$$
\begin{array}{c|l|l|l}
\text{Type of a hyperelliptic surface}&G&m_1,\ldots,m_s&\text{Basis of $\Num(S)$}\\
\hline
1&\Zet_2&2,2,2,2&A/2, B\\
2&\Zet_2\times\Zet_2&2,2,2,2&A/2, B/2\\
3&\Zet_4&2,4,4&A/4, B\\
4&\Zet_4\times\Zet_2&2,4,4&A/4, B/2\\
5&\Zet_3&3,3,3&A/3, B\\
6&\Zet_3\times\Zet_3&3,3,3&A/3, B/3\\
7&\Zet_6&2,3,6&A/6, B
\end{array}
$$
\end{theorem}
Let $\mu=\lcm\{m_1, \ldots, m_s\}$ and let $\gamma=|G|$. Given a hyperelliptic surface, its basis of  $\Num(S)$ consists of divisors $A/\mu$ and $\left(\mu/\gamma\right) B$. We say that $L$ is a line bundle of type $(a,b)$ on a hyperelliptic surface if $L\equiv a\cdot A/\mu+b\cdot(\mu/\gamma) B$.
In  $\Num(S)$  we have that $A^2=0$, $B^2=0$, $AB=\gamma$.

The following proposition holds:
\begin{proposition}[see \cite{Se1990}, Lemma 1.3]\label{Ser2}
Let $D$ be a divisor of type $(a,b)$ on a hyperelliptic surface $S$. 
Then $$D \text{ is ample if and only if } a>0 \text{ and } b>0.$$
\end{proposition}

\vskip 5pt
Now we recall the criterion for a line bundle on a surface to be $k$-very ample, obtained by  M. Beltrametti and A. Sommese in \cite{BeS1988}.

\begin{theorem}[Beltrametti, Sommese]\label{BS}
Let $S$ be a smooth projective surface. Let $L$ be a nef line bundle on $S$ such that $L^2\geq 4k+5$.

Then either  $K_S+L$ is $k$-very ample or there exists an effective divisor 
$D$ satisfying the following conditions:
\begin{enumerate}
\item $L-2D$ is $\Qu$-effective, i.e. there exists an integer  $m>0$ such that $|m(L-2D)|\neq\emptyset.$
\item $D$ contains a subscheme $Z$ of length $k+1$ such that the map
$$H^0(K_S\otimes L) \longrightarrow H^0(K_S\otimes L\otimes \Oo_Z)$$
is not surjective.
\item $LD-k-1\leq D^2<\frac{LD}{2}<k+1$.
\end{enumerate}
\end{theorem}

M. Mella and M. Palleschi in \cite{MP1993} fully characterised $k$-very ampleness of line bundles on hyperelliptic surfaces. For an ample line bundle $L\equiv (a,b)$ they give necessary and sufficient numerical conditions on $a$ and $b$ for each hyperelliptic surface's type.

We will use the sufficient condition for $k$-very ampleness of a line bundle on a  hyperelliptic surface that is implied by \cite{MP1993}, Theorems 3.2-3.4:
\begin{proposition}[Mella, Palleschi]
Let  $S$ be a hyperelliptic surface. Let $L\equiv (a,b)$ be an ample line bundle on $S$. Let $k\in \En$.

If $a\geq k+2$ and $b\geq k+2$ then $L$ is $k$-very ample.
\end{proposition}

In the next section we will prove a condition on the number $r$ for which a pull-back of a $d$-very ample line bundle on a hyperelliptic surface is $k$-very ample on the blow-up of this surface at $r$ very general points.

\section{Main result}

We study  $k$-very ampleness for $k\geq 2$. Case  $k=1$ for a blow-up of a smooth projective surface was considered by M. Coppens, see \cite{Co1995}, Theorem~2. Namely, Coppens proved that on a blow-up of a smooth projective surface at $r$ points in very general position a line bundle  $M=\pi^*(mL)-\sum_{i=1}^r E_i$, where $L$ is an ample line bundle, is $1$-very ample (i.e. very ample) if $m\geq 7$ and $r\leq h^0(mL)-7$. Even if we proved Theorem \ref{k-va} for $k=1$, we would get a weaker result than Coppens.

Our main result is the following

\begin{theorem}\label{k-va}
Let $S$ be a hyperelliptic surface. Let $k\geq 2$, and let $d>(k+1)^2$. Let $L_S\equiv (a,b)$ a line bundle on $S$ with $a\geq d+2$ and $b\geq d+2$.

Let $r\geq 2$.
Let $\pi\colon \widetilde{S}\longrightarrow S$ be the blow-up of $S$ at $r$ points in very general position where
$$r\leq 0.887 \cdot \frac{L_S^2}{(k+1)^2}.$$

Then a line bundle  $L=\pi^*L_S-k\sum_{i=1}^r E_i$ is $k$-very ample on $\widetilde{S}$.
\end{theorem}

Our proof is based on H. Tutaj-Gasi{\'n}ska's ideas from \cite{T-G2005}, Theorem 11. We get a more accurate estimation on the admissible number of points $r$ than in \cite{T-G2005}. This is caused by the fact that for hyperelliptic surfaces we have better estimation of the multi-point Seshadri constants than for arbitrary elliptic quasi-bundle, and on specifics of hyperelliptic surfaces among elliptic fibrations. 

Moreover, assuming that $r\leq c \cdot \frac{L_S^2}{(k+1)^2}$ 
we carefully analysed the conditions for a constant~$c$ to be a maximal possible constant satisfying all conditions imposed by the proof, with any  $\delta>0$. The key restriction for the upper bound of $c$ is given by inequalities \eqref{eq:szac_st_Sesh} and~\eqref{eq:ostatnia}.
The constant $0.887$ is computed to be a round down to the third decimal place of the maximal $c$ satisfying all conditions appearing in the proof.

\begin{proof}
On hyperelliptic surfaces  $K_S\equiv 0$, hence
 $L_S\equiv L_S-K_S\equiv(a,b)$. Obviously, $L_S^2=2ab\geq  2(d+2)^2\geq \left((k+1)^2+3\right)^2$. We prove $k$-very ampleness of $L=\pi^*L_S-k\sum_{i=1}^r E_i$, applying Theorem \ref{BS} to the line bundle 
$$N=L-K_{\widetilde{S}}\equiv \pi^*L_S-(k+1)\sum_{i=1}^r E_i.$$

In the two consecutive lemmas we check that the assumptions of Theorem  \ref{BS} are satisfied, i.e. that  $N^2\geq 4k+5$ and that $N$ is a nef line bundle (we prove that $N$ is in fact ample). Finally, we show that there does not exist an effective divisor $D$ satisfying condition (3) of Theorem \ref{BS}.

\begin{lemma} With the notation above
$$N^2\geq 4k+5.$$
\begin{proof}[Proof of the lemma]
We estimate: \quad 
$ N^2=\left(\pi^*L_S-(k+1)\sum_{i=1}^r
E_i\right)^2=L_S^2-(k+1)^2r\geq
L_S^2-0.887 \cdot \frac{L_S^2}{(k+1)^2}\cdot(k+1)^2=0.113\cdot L_S^2 \geq
0.113\cdot 2\left((k+1)^2+3\right)^2=0.113\cdot 2(k^4+4k^3+12k^2+16k+16)  \geq 0.113\cdot 2\cdot 16(4k+1)\geq 14k+3 \geq 4k+5$.
\end{proof}
\end{lemma}

\begin{lemma}\label{Nszeroka}
$N$ is ample.
\begin{proof}[Proof of the lemma]
By \cite{Fa2015}, Theorem 3.6,  we have
$\varepsilon(L_S, r)\geq\sqrt{\frac{L_S^2}{r}}\sqrt{1-\frac{1}{8r}}.$
We will prove that
$$(\star)\qquad\sqrt{\frac{L_S^2}{r}}\sqrt{1-\frac{1}{8r}}>k+1+\delta$$ where $\delta>0$. Applying an equivalent definition of $r$-point Seshadri constant we will get an assertion of the lemma.

It is enough to show that $(\star)$ holds for the maximal admissible $r$, i.e. for $r= 0.887 \cdot \frac{L_S^2}{(k+1)^2}$. We ask whether
$$\sqrt{\frac{8\cdot 0.887 \cdot \frac{L_S^2}{(k+1)^2}L_S^2- L_S^2}{8\cdot \left(0.887 \cdot \frac{L_S^2}{(k+1)^2}\right)^2}} > k+1+\delta$$
$$\frac {k+1}{0.887}\sqrt{0.887-\frac{(k+1)^2}{8\cdot {L_S^2}}}> k+1+\delta$$

It suffices to check that
\begin{equation}\label{eq:szac_st_Sesh}
(k+1)\left(\frac {1}{0.887}\sqrt{0.887-\frac{(k+1)^2}{8\cdot {2\left((k+1)^2+3\right)^2}}}-1\right) > \delta
\end{equation}

Let $t=k+1$. Computing the derivative of $f(t)=\frac {1}{0.887}\sqrt{0.887-\frac{t^2}{8\cdot {2\left(t^2+3\right)^2}}}$ we see that it is positive, hence $f$ is an increasing function. Evaluating $f$ at the minimal possible $t=3$ (i.e. $k=2$), we get $f(2)\approx 1.0594$. Hence the left hand side of the inequality \eqref{eq:szac_st_Sesh} is an increasing function.
For the minimal  $k=2$ on the left hand side of \eqref{eq:szac_st_Sesh} we get a number slightly bigger than $0.178$ (the difference is on the fourth decimal place). Thus the inequality holds for each $k\geq 2$, if the round down of $\delta$ to the third decimal place is at most $0.178$.

We have proved that  $\varepsilon (L_S,r)>k+1+\delta$ for  $\delta\in(0, \ 0.178]$. Therefore $N$ is ample.
\end{proof}
\end{lemma}

\begin{lemma}There does not exist an effective divisor $D$ such that  $$ND-k-1\leq D^2<\frac{ND}{2}<k+1.$$
\begin{proof}[Proof of the lemma]
Assume that such a divisor exists. Then $D=\pi^*D_S-\sum_{i=1}^r m_iE_i$,
where $m_i:=\mult_{x_i}D_S$. Without loss of generality  $D_S\not \equiv 0$. We consider two cases:
\begin{enumerate}
\item $D^2>0$,
\item $D^2\leq 0$.
\end{enumerate}

Ad. (1). By assumptions of the main theorem
$$r\leq 0.887 \cdot \frac{L_S^2}{(k+1)^2}$$
$$0.113\cdot L_S^2\leq  L_S^2- r \cdot (k+1)^2=N^2$$

Since $N$ is ample, by Hodge Index Theorem $N^2D^2\leq (ND)^2$. Obviously, $N^2\leq N^2D^2$. 
By assumption that a divisor $D$ exists, $\frac{ND}{2}<k+1$.
Moreover, $L_S^2\geq 2((k+1)^2+3)^2$.

Altogether we get
$$0.113\cdot 2((k+1)^2+3)^2\leq 0.113\cdot L_S^2\leq N^2\leq (ND)^2\leq (2k+1)^2.$$

Therefore we have a series of inequalities
$$4k^2+4k+1\geq 0.113\cdot 2(k^4+4k^3+12k^2+16k+16)\geq $$
$$0.226\cdot (4k^2+8k^2+12k^2+16k+16)\geq 0.226\cdot (23k^2+18k+16)>5k^2+4k+3,$$
which gives a contradiction in case $D^2>0$.

\vskip 5pt
Ad. (2). $D^2\leq 0$.

Since $N$ is ample, $ND>0$. Hence $ND\geq 1$. We also have that $ND-k-1\leq D^2$.
Therefore
$$D^2\geq ND-k-1\geq -k.$$

As $D=\pi^*D_S-\sum_{i=1}^r m_iE_i$, we have $D^2=D_S^2-\left(\sum_{i=1}^r m_i\right)^2$. Thus 
$$-k\leq D_S^2-\left(\sum_{i=1}^r m_i\right)^2.$$

Since $ND-k-1\leq D^2$ and $ D^2\leq 0$, we get that $ND\leq k+1$. We compute:  $$\displaystyle ND=\left(\pi^*L_S-(k+1)\sum_{i=1}^r E_i\right).\left(\pi^*D_S-\sum_{i=1}^r m_iE_i\right)=L_SD_S-(k+1)\sum_{i=1}^r m_i.$$

Therefore
$$L_SD_S=ND+(k+1)\sum_{i=1}^r m_i\leq
(k+1)\left(1+\sum_{i=1}^r m_i\right).$$

Since $\left(\sum_{i=1}^r m_i\right)^2\leq D_S^2+k$, we have
$$L_SD_S\leq(k+1)\left(1+\sum_{i=1}^r m_i\right)\leq(k+1)\left(1+D_S^2+k\right).$$

Clearly, $D_S^2\geq 0$. 

If $D_S^2=0$, then  $L_SD_S\leq(k+1)^2$. On the other hand,  $D_S$ is effective and $D_S\not\equiv 0$, hence if $D_S\equiv (\alpha,\beta)$, where $\alpha\geq 0$, $\beta\geq 0$ and $\alpha$ or $\beta$ non-zero, then
$$L_SD_S=a\beta+b\alpha\geq \min\{a,b\}\geq d+2.$$

Therefore
$$(k+1)^2+3\leq d+2\leq L_SD_S\leq(k+1)^2,$$
a contradiction.

\vskip 5pt

If $D_S^2>0$, then by 
$L_SD_S\leq(k+1)\left(1+D_S^2+k\right)$ and Hodge Index Theorem we get
$$L_S^2D_S^2\leq (L_SD_S)^2\leq(k+1)^2\left(1+D_S^2+k\right)^2 \qquad \qquad$$
\hskip 105pt $\|$
$$\qquad \qquad \quad 2abD_S\geq 2(d+2)^2\cdot D_S^2\geq 2((k+1)^2+2)^2\cdot D_S^2\geq 2(k+1)^4\cdot D_S^2.$$

Hence
$$2(k+1)^2\cdot D_S^2\leq\left(k+1+D_S^2\right)^2.$$

We denote $z=D_S^2$, $t=k+1$. We have
$$2t^2 z\leq\left(t+z\right)^2,$$
$$0\leq z^2+\left(2t-2t^2\right)z+t^2,$$
which is a quadratic equation in the variable $z$. Let $z_1(t)=-t+t^2-\sqrt{t^4-2t^3}$, $z_2(t)=-t+t^2+\sqrt{t^4-2t^3}$ be the roots of the equation. We will show that the open interval $(z_1(t), z_2(t))$ contains the closed interval  $\left[1,\frac{1000}{887}t^2\right]$.

We compute the derivative: $z_1'(t)=-1+2t-\frac{3t^2-2t^3}{\sqrt{t^4-2t^3}}$. It is easy to verify that  $z_1'(t)<0$ for all admissible  $t\geq 3$, hence  $z_1(t)$ is a decreasing function. Evaluating $z_1$ at the minimal possible   $t=3$ ($k=2$), we get $z_1(3)\approx 0.804<1$.

Now we compute the derivative: $z_2'(t)-\frac{1000}{887}t^2= -1+2t+\frac{3t^2-2t^3}{\sqrt{t^4-2t^3}}-\frac{1000}{887}t^2$. It is greater than 
 $0$ for all admissible $t\geq 3$, so $z_2(t)-\frac{1000}{887}t^2$ is an increasing function. Evaluating at the minimal possible   $t=3$ ($k=2$), we get the value of approximately $0.001>0$.

Thus we have a contradiction for $0<D_S^2\leq\frac{1000}{887}(k+1)^2$.

\vskip 5pt

Let $D_S^2>\frac{1000}{887}(k+1)^2$. By definition of the multi-point Seshadri constant
$$\varepsilon(L_S,r)\cdot \sum_{i=1}^r m_i\leq L_SD_S.$$

We have already proved that $\varepsilon(L_S,r)\geq k+1+\delta$ so
$$L_SD_S\geq \varepsilon(L_S,r)\cdot \sum_{i=1}^r m_i\geq\left(k+1+\delta\right)\sum_{i=1}^r m_i.$$
On the other hand, we have shown that
$L_SD_S\leq(k+1)\left(1+\sum_{i=1}^r m_i\right)$,  therefore
$$(k+1)\left(1+\sum_{i=1}^r m_i\right)\geq \left(k+1+\delta\right)\sum_{i=1}^r m_i.$$

Setting $t=k+1$, we have
$$t\left(1+\sum_{i=1}^r m_i\right)\geq \left(t+\delta\right)\sum_{i=1}^r m_i.$$
$$\frac{1}{\delta} t\geq \sum_{i=1}^r m_i.$$

Thus
$$L_SD_S\leq t\left(1+\sum_{i=1}^r m_i\right)\leq t\left(1+\frac{1}{\delta}t\right).$$
Squaring both sides we get
$$(L_SD_S)^2\leq t^2\left(1+\frac{1}{\delta}t\right)^2.$$

By Hodge Index Theorem and assumptions $(L_SD_S)^2\geq L_S^2D_S^2> 2\left(t^2+3\right)^2\frac{1000}{887}t^2$, hence 
$$2\left(t^2+3\right)^2\frac{1000}{887}t^2 -t^2\left(1+\frac{1}{\delta}t\right)^2<0,$$
\begin{equation}\label{eq:ostatnia}
2\left(t^2+3\right)^2\frac{1000}{887}-\left(1+\frac{1}{\delta}t\right)^2<0.
\end{equation}

We set the maximal possible by previous computations $\delta= 0.178$ and compute the derivative of $g(t)=2\left(t^2+3\right)^2\frac{1000}{887}-\left(1+\frac{1}{ 0.178}t\right)^2$. We obtain $g'(t)>0$ for $t\geq 3$, hence  $g$ is an increasing function for  $t\geq 3$. Since the value of $g$ for the minimal possible $t=3$ ($k=2$) is positive, we get a contradiction.
\end{proof}
\end{lemma}

We have shown that by Theorem \ref{BS} the divisor $K_{\widetilde{S}}+N$ is  $k$-very ample, but \mbox{$K_{\widetilde{S}}+N=L$.}
\end{proof}

We conclude with a remark.
\begin{remark}
If we improved an estimation of multi-point Seshadri constant of a line bundle on a hyperelliptic surface, then we could easily show that the assertion of main theorem is satisfied with a bigger constant $c$, and therefore the line bundle $L$ is $k$-very ample on the blow-up of a hyperelliptic surface in more very general points. 

However, if we want to apply  Theorem \ref{BS} to $L-K_{\widetilde{S}}$ then a constant $c$, rounded down to the third decimal place, cannot exceed the number  $0.954$, as otherwise the inequality  $N^2\geq 4k+5$ would not hold.
\end{remark}

\subsection*{Acknowledegments} The author would like to thank  Tomasz Szemberg and Halszka Tutaj-Gasi\'nska for many helpful remarks.

\end{document}